\documentclass[fleqn,11pt,a4paper,leqno]{article}

\usepackage{graphicx}
\usepackage{graphics}
\usepackage{latexsym}
\usepackage{amssymb,amsmath}
\usepackage{multicol}
\usepackage{bussproofs}
\usepackage{enumerate}
\usepackage{xcolor}
\usepackage[latin1]{inputenc}
\usepackage{verbatim}

\title{}

\date{}

\numberwithin{equation}{section}

\newtheorem{definicion}{Definition}[section]

\newtheorem{definition}[definicion]{Definition}
\newtheorem{Lemma}[definicion]{Lemma}

\newtheorem{theorem}[definicion]{Theorem}
\newtheorem{Corollary}[definicion]{Corollary}
\newtheorem{corollary}[definicion]{Corollary}
\newtheorem{lema}[definicion]{Lemma}

\newenvironment{Proof}{\noindent\bf Proof \rm}{$\hfill
\square$}

\begin{document}

\title{Semisimple Varieties of Implication Zroupoids}

\author{Juan M. CORNEJO \footnote{I wish to dedicate this work to my little daughter Catalina Cornejo} and Hanamantagouda P. SANKAPPANAVAR}

\numberwithin{equation}{section}

\maketitle

\begin{abstract}
It is a well known fact that Boolean algebras can be defined using only implication and a constant.   In 2012, this result was extended to De Morgan algebras in \cite{sankappanavarMorgan2012} which led Sankappanavar to introduce, and investigate, the variety $\mathbf{I}$ of implication zroupoids generalizing De Morgan algebras.   His investigations were continued in \cite{CoSa2015aI} and \cite{CoSa2015a} in which several new subvarieties of $\mathbf{I}$ were introduced and their relationships with each other and with the varieties of  \cite{sankappanavarMorgan2012} were explored.

The present paper is a continuation of \cite{sankappanavarMorgan2012} and \cite{CoSa2015aI}. 
The main purpose of this paper is to determine the simple algebras in $\mathbf{I}$.  It is shown that there are exactly five simple algebras in $\mathbf{I}$. 
From this description we deduce that the semisimple subvarieties of $\mathbf{I}$ are precisely the subvarieties of the variety
generated by these 5 simple I-zroupoids and are locally finite.  It also follows that the lattice of semisimple subvarieties of $\mathbf{I}$ is isomorphic to 
the direct product of a 4-element Boolean lattice and a 4-element chain.
\end{abstract}

\thispagestyle{empty}

\section{Introduction}

It is a well known fact that Boolean algebras can be defined using only implication and a constant.   In 2012, this result was extended to De Morgan algebras in \cite{sankappanavarMorgan2012} which led Sankappanavar to introduce, and investigate, the variety $\mathbf{I}$ of implication zroupoids generalizing De Morgan algebras.  His investigations were continued in \cite{CoSa2015aI} and \cite{CoSa2015a} in which several new subvarieties of $\mathbf{I}$ were introduced and their relationships with each other and with the varieties of  \cite{sankappanavarMorgan2012} were explored.
   
The present paper is a continuation of \cite{sankappanavarMorgan2012} and \cite{CoSa2015aI}.
The main purpose of this paper is to determine the simple algebras  in $\mathbf{I}$.  It is shown that there are five simple algebras in $\mathbf{I}$, namely:
\begin{itemize}
\item
    the 2-element trivial implication zroupoid $\mathbf{2_z}$, where $x \to y :=0$
 \item
       the 2-element $\lor$-semilattice $\mathbf{2_s}$ with the least element $0$,
\item
      the 2-element Boolean algebra $\mathbf{2_b}$,
\item
     the 3-element Kleene algebra $\mathbf{3_k}$, and
\item
     the 4-element De Morgan algebra $\mathbf{4_d}$.
\end{itemize} 
From this description it follows that the semisimple subvarieties of $\mathbf{I}$ are precisely the subvarieties of the variety
$\mathbb{V(\mathbf{2_z}, \mathbf{2_s}, \mathbf{2_b}, \mathbf{3_k}, \mathbf{4_d})}$ generated by the 5 simple I-zroupoids and are locally finite.  It also follows that the lattice of semisimple varieties of implication zroupoids is isomorphic to 
the direct product of a 4-element Boolean lattice and a 4-element chain.

The method employed in this paper to describe simple algebras of $\mathbf{I}$ will, we hope,  inspire a way for solving the problem of characterizing all the subdirectly irreducible members of $\mathbf{I}$ and, hence, for obtaining information about the lattice of subvarieties of the variety of implication zroupoids.

We include an appendix with some of the longer proofs 
so that the reader can focus on the main ideas and results of the paper.

\section{Preliminaries} \label{SA}

For the concepts not defined, and results not proved, in this paper we refer the reader to the textbooks \cite{balbesDistributive1974},
\cite{burrisCourse1981}, and \cite{Ra74}.  The familiarity with \cite{CoSa2015aI} and \cite{sankappanavarMorgan2012}, although not necessary, is helpful to the reader.
	We recall from \cite{sankappanavarMorgan2012} the following definition which is central to this paper.
	
\begin{definition}
An algebra $\mathbf A = \langle A, \to, 0 \rangle$, where $\to$ is binary and $0$ is a constant, is called a {\it zroupoid}.
A zroupoid $\mathbf A = \langle A, \to, 0 \rangle$ is an {\it Implication zroupoid} {\rm(I}-zroupoid, for short{\rm)} if $\mathbf A$ satisfies:
\begin{itemize}
	\item[\rm{(I)}] 	$(x \to y) \to z \approx [(z' \to x) \to (y \to z)']'$, where $x' : = x \to 0$
	\item[{\rm (I$_{0}$)}]  $ 0'' \approx 0$.
\end{itemize}

Throughout this paper $\mathbf{I}$ denotes the variety of implication zroupoids.
\end{definition}
In this paper we use the characterizations of De Morgan algebras, Kleene algebras and Boolean algebras obtained in \cite{sankappanavarMorgan2012} as definitions.

\begin{definition}
An implication zroupoid $\mathbf A = \langle A, \to, 0 \rangle$ is a {\it De Morgan algebra} {\rm(}$\mathbf{DM}$-algebra for short{\rm)} if $\mathbf A$ satisfies the axiom:
\begin{itemize}
	\item
	[{\rm(DM)}]  $(x \to y) \to x \approx x. $
\end{itemize}

A $\mathbf{DM}$-algebra $\mathbf A = \langle A, \to, 0 \rangle$ is a {\it Kleene algebra} {\rm(}$\mathbf{KL}$-algebra for short{\rm)} if
$\mathbf A$ satisfies the axiom:
\begin{itemize}
	\item[{\rm(KL$_1$)}] $(x \to x) \to (y \to y)' \approx x \to x$
\end{itemize}
or, equivalently,
\begin{itemize}
	\item[{\rm(KL$_2$)}] $(y \to y) \to (x \to x) \approx x \to x$.
\end{itemize}

A $\mathbf{DM}$-algebra $\mathbf A = \langle A, \to, 0 \rangle$ is a {\it Boolean algebra} {\rm(}$\mathbf{BA}$-algebra for short{\rm)} if
$\mathbf A$ satisfies the axiom:
\begin{itemize}
	\item[{\rm(BA)}] $x \to x \approx 0'$.
\end{itemize}
We denote by $\mathbf{DM}$, $\mathbf{KL}$, $\mathbf{BA}$ respectively the varieties of $\mathbf{DM}$-algebras,
$\mathbf{KL}$-algebras and $\mathbf{BA}$-algebras. 
\end{definition}


\medskip

\section{Some properties of Implication Zroupoids in $\mathbf{I_{2,0}}$}

Recall from \cite{sankappanavarMorgan2012} that $\mathbf{I_{2,0}}$ denotes the subvariety of  $\mathbf{I}$
defined by the identity: $x'' \approx x$.  In this section, in addition to recalling some known results, we present  some new properties of $\mathbf{I_{2,0}}$ which, being of interest in their own right, will be useful later.

\medskip

\begin{lema} {\rm \cite{sankappanavarMorgan2012}} \label{general_properties}
	Let $\mathbf A \in \mathbf{I_{2,0}}$. Then
	\begin{enumerate}[{\rm (a)}]
		\item $x' \to 0' \approx 0 \to x$ \label{cuasiConmutativeOfImplic2}
		\item $0 \to x' \approx x \to 0'$. \label{cuasiConmutativeOfImplic}
	\end{enumerate}
\end{lema}

\begin{lema} {\rm \cite[Theorem 8.15]{sankappanavarMorgan2012}} \label{general_properties_equiv}
	Let $\mathbf A$ be an {\bf{I}}-zroupoid. Then the following are equivalent:
	\begin{enumerate}[{\rm (a)}]
		\item $0' \to x \approx x$ \label{TXX} 
		\item $x'' \approx x$
		\item $(x \to x')' \approx x$ \label{reflexivity}
		\item $x' \to x \approx x$. \label{LeftImplicationwithtilde}
	\end{enumerate}
\end{lema}

The following lemma, whose proof is given in the Appendix, plays a crucial role in the rest of the paper.

\begin{lema} \label{general_properties2}
	Let $\mathbf A \in \mathbf{I_{2,0}}$. Then:
	\begin{enumerate}[{\rm (1)}]
		\item $(x \to 0') \to y \approx (x \to y') \to y$ \label{281014_05} 
		\item $(0 \to x') \to (y \to x) \approx y \to x$ \label{281014_07} 
		\item $x \to (0 \to x)' \approx x'$ \label{291014_02} 
		\item $(y \to x)' \approx (0 \to x) \to (y \to x)'$ \label{291014_06} 
		\item $[x \to (y \to x)']' \approx (x \to y) \to x$ \label{291014_09} 
		\item $(y \to x) \to y \approx (0 \to x) \to y$ \label{291014_10} 
		\item $0 \to x \approx 0 \to (0 \to x)$ \label{311014_03} 
		\item $0 \to (x \to y')' \approx 0 \to (x' \to y)$ \label{191114_05} 
		\item $[(x \to y) \to x] \to [(y \to x) \to y] \approx x \to y$ \label{271114_03} 
		\item $[x' \to (0 \to y)]' \approx (0 \to x) \to (0 \to y)'$ \label{031114_06} 
		\item $0 \to (0 \to x)' \approx 0 \to x'$ \label{031114_07} 
		\item $0 \to (x' \to y)' \approx x \to (0 \to y')$ \label{071114_01} 
		\item $0 \to (x \to y) \approx x \to (0 \to y)$ \label{071114_04} 
		\item $[(0 \to x) \to y] \to x \approx [(y \to x) \to (0 \to x)']'$ \label{291014_03} 
		\item $(x \to y) \to (0 \to y)' \approx (x \to y)'$ \label{291014_07} 
		\item $x \to [(y \to z') \to x]' \approx (0 \to y) \to \{(0 \to z) \to x'\}$ \label{181114_10} 
		\item $(x \to y) \to y' \approx y \to (x \to y)'$ \label{071114_05} 
		\item $0 \to \{(0 \to x) \to y'\} \approx x \to (0 \to y')$ \label{071114_03} 
		\item $[(x \to 0') \to y]' \approx (0 \to x) \to y'$ \label{071114_02} 

		\item $[(0 \to x) \to y] \to x \approx y \to x$ \label{291014_08}
		\item $(0 \to x) \to (0 \to y) \approx x \to (0 \to y)$ \label{311014_06}
		\item $0 \to [(x \to y) \to z']' \approx (x \to y)' \to (0 \to z)$ \label{031114_03}
		\item $[(x' \to y') \to (y \to x)'] \to (y \to x) \approx y \to x$ \label{311014_01}
		\item $[0 \to (x' \to y')'] \to (y \to x) \approx y \to x$ \label{311014_02}
		\item $[0 \to (x \to y)] \to x \approx (0 \to y) \to x$ \label{311014_04}
		\item $(0 \to x) \to (x \to y) \approx x \to (x \to y)$ \label{311014_05}
		\item $(0 \to x) \to [y \to (0 \to x)]' \approx [y \to (0 \to x)]'$ \label{311014_07}
		\item $[(x \to y) \to (0 \to x)]' \approx (0 \to x) \to [y \to (0 \to x)']$ \label{311014_08}
		\item $(x \to y) \to (0 \to x) \approx y \to (0 \to x)$ \label{031114_01}
		\item $y \to (y \to x) \approx [x \to (0 \to y)] \to (y \to x)$ \label{031114_02}
		\item $x \to y \approx x \to (x \to y)$ \label{031114_04}
		\item $[\{x \to (0 \to y)\} \to z]' \approx z \to [(x \to y) \to z]'$ \label{171114_01} 
		
		\item $[0 \to (x \to y)'] \to y \approx (0 \to x) \to y$ \label{181114_01}
		
		\item $[\{x \to (0 \to (y \to z)')\} \to z]' \approx (z' \to x) \to [z \to (y \to z)']$ \label{181114_02} 
		\item $x \to [\{y \to (z \to x)'\} \to x]' \approx (x' \to y) \to [x \to (z \to x)']$ \label{181114_03}
		\item $[0 \to (x \to y)] \to y' \approx y \to (x \to y)'$ \label{181114_11}
		\item $(x' \to y) \to [x \to (z \to x)'] \approx (0 \to y) \to [x \to (z \to x)']$ \label{181114_12}  
		\item $[x' \to (y \to z)] \to [x \to (y \to x)'] \approx x \to [(z \to y) \to x]'$ \label{181114_13}
		\item $x \to [(y \to z) \to x]' \approx [0 \to (z \to y)] \to [x \to (z \to x)']$ \label{181114_14}
		\item $[0 \to (x \to y)] \to [z \to (x \to z)'] \approx (0 \to y) \to [z \to (x \to z)']$ \label{181114_15}
		\item $x \to [(y \to z) \to x]' \approx (0 \to y) \to [x \to (z \to x)']$ \label{181114_16}
		\item $0 \to [(0 \to x) \to y] \approx x \to (0 \to y)$ \label{191114_02}
		\item $x \to (y \to x)' \approx (y \to 0') \to x'$ \label{031214_16}

		\item if $(x \to y') \to x = x$ then $(x' \to y) \to x' = x'$. \label{130315_02}    
		\item $[(0 \to x')\to y] \to [0 \to (z \to x)'] \approx [\{y \to (0 \to z)\} \to (x \to 0')']'$ \label{130315_08}
		\item $[\{x \to (0 \to y)\} \to z]' \approx (0 \to x) \to [z \to (y \to z)']$ \label{130315_09}
		\item $[(0 \to x') \to y] \to [0 \to (z \to x)'] \approx y \to (0 \to (z \to x)')$ \label{130315_10}
		\item $[(0 \to x) \to y] \to (0 \to z) \approx 0 \to [(x \to y) \to z]$ \label{130315_11}
		\item $(x' \to y) \to [0 \to (z \to x)'] \approx y \to [0 \to (z \to x)']$ \label{130315_12}  
		\item $y \to [0 \to (z \to x)'] \approx 0 \to [(x' \to y) \to (z \to x)']$ \label{130315_13}
		\item $(0 \to x) \to  [\{0 \to (y \to z)\} \to u'] \approx [0 \to \{(x \to y) \to z\}] \to u'$ \label{140315_01}
		\item $[0 \to \{(x \to y) \to z\}]' \approx (0 \to x) \to [(0 \to y') \to (0 \to z)']$ \label{140315_02}  
		\item $0 \to [(x \to y) \to (z \to x')']\ \approx 0 \to [y \to (z' \to x)]$ \label{140315_03}
		\item $0 \to [(x \to y)' \to z] \approx 0 \to [x \to (y' \to z)]$ \label{140315_04}
		\item $[\{(x \to y) \to z\} \to \{0 \to (y \to z)\}'] \to [(x \to y) \to z] \approx (x \to y) \to z $ \label{140315_05}
		\item $[0 \to (x \to y)'] \to (z \to y) \approx (0 \to x) \to (z \to y)$ \label{170315_01}
		\item $[x \to (y' \to z)'] \to x \approx [y \to \{0 \to (0 \to z)'\}] \to x$ \label{250315_01}
		\item $[(0 \to x) \to y] \to (z \to x) \approx y \to (z \to x)$ \label{250315_02} 
		\item $x \to [(y \to x) \to y] \approx x \to y$ \label{250315_03}
		\item $[(x \to y) \to (y \to z)]' \approx (0 \to x) \to (y \to z)'$  \label{250315_04}
		\item $(x \to y) \to (y \to x) \approx y \to x$  \label{250315_05}
		\item $[\{(x \to y) \to z\} \to (z' \to x)'] \to [(x \to y) \to z] \approx (x \to y) \to z $  \label{250315_06} 
		\item $[\{(x \to y) \to z\} \to \{z' \to (y \to x)\}'] \to [(x \to y) \to z] \approx (x \to y) \to z$.  \label{250315_07}
		
	\end{enumerate}
\end{lema}


\section{Simple algebras in the variety $\mathbf{I_{2,0}}$}

In this section we will prove that if $\mathbf A \in \mathbf{I_{2,0}}$ is a simple algebra with $|A| \geq 3$, then $\mathbf{A} \in \mathbf{DM}$.

\medskip

\begin{definition} \label{definition_relat_R1}
	Let $\mathbf A \in \mathbf{I_{2,0}}$. We define the relation $R_1$ on $\mathbf A$ as follows:
	$$x R_1 y \mbox{ if and only if } (x \to y') \to x = x \mbox{ and } (y \to x') \to y = y.$$
\end{definition}

\begin{Lemma} \label{lemma_equivalence_relation_simpleAlg}
	Let $\mathbf A \in \mathbf{I_{2,0}}$. Then $R_1$ is a equivalence relation on $\mathbf A$.
\end{Lemma}

\begin{Proof}
	Let $a \in A$. By Lemma \ref{general_properties_equiv} (\ref{LeftImplicationwithtilde}) we have that $(a \to a') \to a = (a'' \to a') \to a = a' \to a = a$. Then the relation $R_1$ is reflexive. The symmetry is immediate.
Now we show that $R_1$ is transitive.  Let $a,b,c \in A$ and assume that $a R_1 b$ and $b R_1 c$. This implies that
	 \begin{equation} \label{130315_01}
	(a \to b') \to a = a,
	\end{equation}
	\begin{equation} \label{130315_06}
	(b \to a') \to b = b,
	\end{equation}
	\begin{equation} \label{130315_03}
	(b \to c') \to b = b
	\end{equation}
	and
	\begin{equation} \label{130315_07}
	(c \to b') \to c = c.
	\end{equation}
	
By Lemma \ref{general_properties2} (\ref{130315_02}) and Lemma (\ref{130315_01}) we have 
	\begin{equation} \label{130315_05}
	a' = (a' \to b) \to a'.
	\end{equation}
	Therefore,
	$$
	\begin{array}{lcll}
	a' & = & (a' \to b) \to a' & \mbox{by (\ref{130315_05})} \\
	& = & (0 \to b) \to a'  & \mbox{by Lemma \ref{general_properties2} (\ref{291014_10})} \\
	& = & (0 \to ((b \to c') \to b)) \to a' & \mbox{by (\ref{130315_03})} \\
	& = & [0 \to \{(0 \to c') \to b\}]  \to a' & \mbox{by Lemma \ref{general_properties2} (\ref{291014_10})} \\
	& = & [(0 \to c') \to (0 \to b)] \to a' & \mbox{by Lemma \ref{general_properties2} (\ref{071114_04})} \\
	& = & [\{a'' \to (0 \to c')\} \to \{(0 \to b) \to a'\}']' & \mbox{by (I)} \\
	& = & [\{a \to (0 \to c'\}] \to \{(0 \to b) \to a'\}']'  & \mbox{} \\
	& = & [\{a \to (0 \to c')\} \to \{(a' \to b) \to a'\}']' & \mbox{by Lemma \ref{general_properties2} (\ref{291014_10})} \\
	& = & [\{a \to (0 \to c')\} \to a'']' & \mbox{by (\ref{130315_05})} \\
	& = & [\{a \to (0 \to c')\} \to a]' & \mbox{} \\
	& = & [\{0 \to (0 \to c')\} \to a]' & \mbox{by Lemma \ref{general_properties2} (\ref{291014_10})} \\
	& = & [(0 \to c') \to a]' & \mbox{by Lemma \ref{general_properties2} (\ref{311014_03})} \\
	& = & [(a \to c') \to a]' & \mbox{by Lemma \ref{general_properties2} (\ref{291014_10})}
	\end{array}
	$$
	and, hence we get  $$a = a'' = [(a \to c') \to a]'' = (a \to c') \to a.$$
	Similary, from (\ref{130315_06}) and (\ref{130315_07}), we conclude that $$(c \to a') \to c = c,$$
	implying $a R_1 c$.	
\end{Proof}

\begin{Lemma} \label{lemma_congruence_relation_simpleAlg}
	Let $\mathbf A \in \mathbf{I_{2,0}}$. Then $R_1$ is a congruence.
\end{Lemma}

\begin{Proof} By Lemma \ref{lemma_equivalence_relation_simpleAlg}, $R_1$ is an equivalence relation.
	Let $a,b,c,d \in A$ satisfying $a R_1 b$ and $c R_1 d$. Then
	\begin{equation}\label{160315_01}
	(a \to b') \to a = a,
	\end{equation}
	\begin{equation} \label{170315_02}
	(b \to a') \to b = b,
	\end{equation}
	\begin{equation} \label{250315_08}
	(c \to d') \to c = c
	\end{equation}
	and
	\begin{equation} \label{250315_12}
	(d \to c') \to d = d.
	\end{equation}
	From (\ref{160315_01}) and (\ref{250315_08}), and Lemma \ref{general_properties2} (\ref{130315_02}), we have that
	\begin{equation}\label{160315_02}
	(a' \to b) \to a' = a'.
	\end{equation}
	and
	\begin{equation}\label{250315_09}
	(c' \to d) \to c' = c'.
	\end{equation}
	Consequently,
	$$
	\begin{array}{lcll}
	a' \to [(b \to a')'] & = &[a' \to [(b \to a')']]''  & \mbox{} \\
	& = & [(a' \to b) \to a']' & \mbox{by Lemma \ref{general_properties2} (\ref{291014_09})} \\
	& = & a'' & \mbox{using (\ref{160315_02})} \\
	& = & a. & \mbox{}
	\end{array}
	$$
	Therefore,
	\begin{equation}\label{160315_03}
	a' \to [(b \to a')'] = a.
	\end{equation}
	Hence we have
	$$
	\begin{array}{lcll}
	(0 \to b) \to (a \to c) & = & (b' \to 0') \to (a \to c) & \mbox{by Lemma \ref{general_properties} (\ref{cuasiConmutativeOfImplic})
	} \\
	& = & [(b \to 0) \to 0'] \to (a \to c) & \mbox{} \\
	& = & [(b \to 0'') \to 0'] \to (a \to c) & \mbox{} \\
	& = & [(b \to 0') \to 0'] \to (a \to c) & \mbox{by Lemma \ref{general_properties2} (\ref{281014_05})} \\
	& = & [(b \to 0') \to (a \to c)'] \to (a \to c) & \mbox{by Lemma \ref{general_properties2} (\ref{281014_05})} \\
	& = & [\{b \to (a \to c)''\} \to (a \to c)'] \to (a \to c) & \mbox{by Lemma \ref{general_properties2} (\ref{281014_05})} \\
	& = & [\{b \to (a \to c)\} \to (a \to c)'] \to (a \to c) & \mbox{} \\
	& = & [\{b \to (a \to c)\} \to 0'] \to (a \to c) & \mbox{by Lemma \ref{general_properties2} (\ref{281014_05})} \\
	& = & [\{b \to (a \to c)\} \to (0 \to 0)] \to (a \to c) & 
	\\
	& = & [\{0 \to (b \to (a \to c)\}] \to (0 \to 0)] \to (a \to c) & \mbox{by Lemma \ref{general_properties2} (\ref{311014_06})} \\
	& = & [0 \to [[0 \to (b \to (a \to c))] \to 0]] \to (a \to c) & \mbox{by Lemma \ref{general_properties2} (\ref{071114_04})} \\
	& = & [(a \to c) \to [[0 \to (b \to (a \to c))] \to 0]] \to (a \to c) & \mbox{by Lemma \ref{general_properties2} (\ref{291014_10})} \\
	& = & [(a \to c) \to [0 \to (b \to (a \to c))]'] \to (a \to c) & \mbox{} \\
	& = & [(a \to c) \to [0 \to [(b \to a')' \to c]]'] \to (a \to c) & \mbox{by Lemma \ref{general_properties2} (\ref{140315_04}) with } \\
	&  &  & x = b, y = a', z = c
	\end{array}
	$$
	Then
	\begin{equation} \label{160315_04}
	(0 \to b) \to (a \to c) = [(a \to c) \to [0 \to [(b \to a')' \to c]]'] \to (a \to c).
	\end{equation}
	By (\ref{160315_04}) and (\ref{160315_03}), we have
	\begin{equation} \label{160315_05}
	(0 \to b) \to (a \to c) = [([a' \to [(b \to a')']] \to c) \to [0 \to [(b \to a')' \to c]]'] \to ([a' \to [(b \to a')']] \to c).
	\end{equation}
    From (\ref{160315_05}) and Lemma \ref{general_properties2} (\ref{140315_05}) with $x = a', y = (b \to a')'$ and $z = c$ we conclude that 
    $$(0 \to b) \to (a \to c) = [a' \to [(b \to a')']] \to c. $$
	Therefore, by (\ref{160315_03}), we have
	\begin{equation}
	(0 \to b) \to (a \to c) = a \to c.
	\end{equation}
	Consequently, in view of Lemma \ref{general_properties2} (\ref{170315_01}), we can conclude that $a \to c = [0 \to (b \to c)'] \to (a \to c)$. From Lemma \ref{general_properties2} (\ref{291014_10}), we obtain
	$$
	a \to c = [(a \to c) \to (b \to c)'] \to (a \to c).
	$$
	Similarly, in view of (\ref{170315_02}), we have
	$$
	b \to c = [(b \to c) \to (a \to c)'] \to (b \to c).
	$$
	Thus we conclude 
	\begin{equation} \label{250315_14}
	a \to c\  R_1 \ b \to c.
	\end{equation}

	Moreover, from Lemma \ref{general_properties2}  (\ref{281014_07}), we have 
	\begin{equation} \label{250315_10}
	b \to c = (0 \to c') \to (b \to c).
	\end{equation}
	Now, we observe that
	$$
	\begin{array}{lcll}
	b \to c & = & (0 \to c') \to (b \to c) & \mbox{by (\ref{250315_10})} \\
	& = & [(0 \to c') \to (b \to c)]' \to [(0 \to c') \to (b \to c)] & \mbox{by Lemma \ref{general_properties_equiv} (\ref{LeftImplicationwithtilde})} \\
	& = & [(0 \to ((c' \to d) \to c')) \to (b \to c)]' \to [(0 \to c') \to (b \to c)] & \mbox{by (\ref{250315_09})} \\
	& = & [(0 \to ((0 \to d) \to c')) \to (b \to c)]' \to [(0 \to c') \to (b \to c)] & \mbox{by Lemma \ref{general_properties2}  (\ref{291014_10})} \\
	& = & [((0 \to d) \to (0 \to c')) \to (b \to c)]' \to [(0 \to c') \to (b \to c)] & \mbox{by Lemma \ref{general_properties2}  (\ref{071114_04})} \\
	& = & [(d \to (0 \to c')) \to (b \to c)]' \to [(0 \to c') \to (b \to c)]  & \mbox{by Lemma \ref{general_properties2}  (\ref{311014_06})} \\
	& = & [[(b \to c)' \to d] \to [(0 \to c') \to (b \to c)]']'' \to [(0 \to c') \to (b \to c)] & \mbox{by (I)} \\
	& = & [[(b \to c)' \to d] \to [(0 \to c') \to (b \to c)]'] \to [(0 \to c') \to (b \to c)] & \mbox{} \\
	& = & [[(b \to c)' \to d] \to 0'] \to [(0 \to c') \to (b \to c)] & \mbox{by  Lemma \ref{general_properties2} (\ref{281014_05})} \\
	& = & [0 \to [(b \to c)' \to d]'] \to [(0 \to c') \to (b \to c)] & \mbox{by Lemma \ref{general_properties} (\ref{cuasiConmutativeOfImplic})} \\
	& = & [[(0 \to c') \to (b \to c)] \to [(b \to c)' \to d]'] \to [(0 \to c') \to (b \to c)]  & \mbox{by Lemma \ref{general_properties2}  (\ref{291014_10})} \\
	& = & [[b \to c] \to [(b \to c)' \to d]'] \to [b \to c] & \mbox{by (\ref{250315_10})} \\
	& = & [[b \to c] \to [0 \to (0 \to d)']] \to [b \to c]  & \mbox{by Lemma \ref{general_properties2}  (\ref{250315_01}) with } \\
	&  &   &  x = b \to c, y = b \to c, z = d \\
	& = & [0 \to [0 \to (0 \to d)']] \to [b \to c] & \mbox{by Lemma \ref{general_properties2}  (\ref{291014_10})} \\
	& = & [0 \to (0 \to d)'] \to [b \to c]  & \mbox{by Lemma \ref{general_properties2}  (\ref{031114_04})} \\
	& = & (0 \to d')  \to [b \to c] & \mbox{by Lemma \ref{general_properties2}  (\ref{031114_07})} \\
	& = & [(b \to c)'' \to [d' \to (b \to c)]']' & \mbox{by (I)} \\
	& = & [(b \to c) \to [d' \to (b \to c)]']' & \mbox{} \\
	& = & [(0 \to d) \to [(b \to c) \to (0 \to (b \to c))']]' & \mbox{by Lemma \ref{general_properties2}  (\ref{181114_16}) with } \\
	&  &  &  x = b \to c, y = d, z = 0 \\
	& = & [(0 \to d) \to (b \to c)']' & \mbox{by Lemma \ref{general_properties2}  (\ref{250315_03}) with } \\
	&  &  & x = b \to c, y = 0 \\
	& = & [(d \to b) \to (b \to c)]'' & \mbox{by  Lemma \ref{general_properties2} (\ref{250315_04})} \\
	& = & (d \to b) \to (b \to c) & \mbox{}
	\end{array}
	$$
	Consequently,
	\begin{equation} \label{250315_11}
	b \to c = (d \to b) \to (b \to c).
	\end{equation}
	Hence,
	$$
	\begin{array}{lcll}
	b \to c & = & (d \to b) \to (b \to c) & \mbox{by (\ref{250315_11})} \\
	& = & [[(d \to b) \to (b \to c)] \to [(b \to c)' \to (b \to d)]'] \to [(d \to b) \to (b \to c)] & \mbox{by Lemma \ref{general_properties2}  (\ref{250315_07}) with }  \\
	&  &  & x = d, y = b, z = b \to c \\
	& = & [[b \to c] \to [(b \to c)' \to (b \to d)]'] \to [b \to c] & \mbox{by (\ref{250315_11})} \\
	& = & [0 \to [(b \to c)' \to (b \to d)]'] \to (b \to c) & \mbox{by  Lemma \ref{general_properties2} (\ref{291014_10})} \\
	& = & [0 \to [[(b \to d)' \to (b \to c)] \to [0 \to (b \to d)]']''] \to (b \to c) & \mbox{by (I)} \\
	& = & [0 \to [[(b \to d)' \to (b \to c)] \to [0 \to (b \to d)]']] & \mbox{} \\
	& = & [(b \to c) \to [0 \to (0 \to (b \to d))']] \to (b \to c) & \mbox{by Lemma \ref{general_properties2}  (\ref{130315_13}) with }  \\
	&  &  & x = b \to d, y = b \to c, z = 0 \\
	& = & [(b \to c) \to [0 \to(b \to d)']] \to (b \to c) & \mbox{by  Lemma \ref{general_properties2} (\ref{031114_07})} \\
	& = & [0 \to [0 \to(b \to d)']] \to (b \to c) & \mbox{by  Lemma \ref{general_properties2} (\ref{291014_10})} \\
	& = & (0 \to(b \to d)') \to (b \to c) & \mbox{by  Lemma \ref{general_properties2} (\ref{031114_04})} \\
	& = & ((b \to c) \to(b \to d)') \to (b \to c) & \mbox{by  Lemma \ref{general_properties2} (\ref{291014_10})}
	\end{array}
	$$
	In a similar way, by (\ref{250315_12}), we have that
	$$b \to d = ((b \to d) \to(b \to c)') \to (b \to d).$$
	Then
	\begin{equation} \label{250315_13}
	b \to c \  R_1\   b \to d.
	\end{equation}
	In view of (\ref{250315_14}), (\ref{250315_13}) and Lemma (\ref{lemma_equivalence_relation_simpleAlg}), we conclude $a \to c \  R_1 \  b \to d.$
\end{Proof}

\begin{Lemma}  \label{250315_16}
	Let $\mathbf A \in \mathbf{I_{2,0}}$ be a simple algebra with $|A| \geq 3$, and let $a,b \in A$ with $a,b \not= 0$. Then
	\begin{enumerate}[{\rm(a)}]
		\item If $0 = 0'$, then $0 \to a \not=0$, $a \to 0 \not=0$, $a \to b \not= 0$. \label{250315_15}
		\item $0' \not= 0.$ 
		             \label{250315_17}
	\end{enumerate}
\end{Lemma}

\begin{Proof}
	\begin{enumerate}[(a)]
		\item Let $a,b \in A$ with $a,b \not= 0$. Assume that $0 \to a = 0$. Then $0 = 0 \to a = 0' \to a = a$ by Lemma \ref{general_properties_equiv} (\ref{TXX}), which is a contradiction, implying $0 \to a \not= 0$.
		
		If $a \to 0 = 0$ then $a = a'' = (a \to 0)' = 0' = 0$, a contradiction.
		
		If $a \to b = 0$, then, by Lemma \ref{general_properties2} (\ref{031114_04}), we have $0 = a \to b = a \to (a \to b) = a \to 0$, a contradiction.
		
		\item Assume that $0 = 0'$. Considerer the equivalence relation, $R_2$ associated to the partition $\{\{0\},A^\ast\}$ of $A$ with $A^\ast = \{c \in A:\ c \not= 0\}$.
		Since $|A| \geq 3$, we observe that $R_2 \not= \Delta, A \times A$.  Now, let $a,b,c,d \in A$ satisfying $a R_2 b$ and  $c R_2 d$. From $0' = 0$ and item (\ref{250315_15}) we conclude that $R_2$ is a congruence, which is a contradiction, because $\mathbf A$ is simple. Thus $0' \not= 0$.
	\end{enumerate}
\end{Proof}

We are now ready to present our main theorem of this section.

\begin{theorem} \label{simpleAlgebras_theorem1}
	Let $\mathbf A \in \mathbf{I_{2,0}}$ be a simple algebra with $|A| \geq 3$. Then $\mathbf A \in \mathbf{DM}$.
\end{theorem}

\begin{Proof}
	We know, in view of Lemma \ref{lemma_congruence_relation_simpleAlg}, that the relation $R_1,$ defined in Definition \ref{definition_relat_R1}, is a congruence on $\mathbf A$.  We also know that $0' \not= 0$ from Lemma \ref{250315_16} (\ref{250315_17}). Observe that
	$$(0 \to 0'') \to 0 = (0 \to 0) \to 0 = 0'' = 0$$
	and
	$$(0' \to 0') \to 0' = 0' \to 0' = 0,'$$
	by Lemma \ref{general_properties_equiv} (\ref{TXX}).
	Hence $0 R_1 0'$. Since $\mathbf A$ is simple and $0' \not= 0$, it follows that $R_1 = A \times A$. Therefore,
	$$A \models (x \to y') \to x \approx x.$$
	Thus, letting $a,b \in A$, we can conclude that
	$(a \to b) \to a = (a \to b'') \to a = a$. Consequently, we have proved
	$$A \models (x \to y') \to x \approx x,$$ implying that
	$\mathbf A \in \mathbf{DM}$.
\end{Proof}

\section{Semisimple Implication Zroupoids}\par

In this section we dermine the semisimple subvarieties of $\mathbf{I}$.

Recall from \cite{sankappanavarMorgan2012} that there are three $2$-element algebras in $\mathbf{I}$, namely $\mathbf{2_z}$, $\mathbf{2_s}$, $\mathbf{2_b}$, whose $\to$
operations are, respectively, as follows:   \\

\begin{minipage}{0.3 \textwidth}
	\begin{tabular}{r|rr}
		$\to$: & 0 & 1\\
		\hline
		0 & 0 & 0 \\
		1 & 0 & 0
	\end{tabular} \hspace{.5cm}	
	
\end{minipage}
\begin{minipage}{0.3 \textwidth}
	\begin{tabular}{r|rr}
		$\to$: & 0 & 1\\
		\hline
		0 & 0 & 1 \\
		1 & 1 & 1
	\end{tabular} \hspace{.5cm}	
	
\end{minipage}
\begin{minipage}{0.3 \textwidth}
	\begin{tabular}{r|rr}
		$\to$: & 0 & 1\\
		\hline
		0 & 1 & 1 \\
		1 & 0 & 1
	\end{tabular} \hspace{.5cm}
	
\end{minipage}

\vskip 1cm

It is clear that $\mathbf{2_z}$, $\mathbf{2_s}$, $\mathbf{2_b}$ are the only $2$-element simple algebras  in $\mathbf{I}$.

\medskip

Recall, also, that the $3$-element (Kleene) algebra, namely $\mathbf 3_k$, and the $4$-element algebra, namely $\mathbf 4_d$, are in $\mathbf{DM}$.  Their $\to$ operations are, respectively, as follows: \\

\begin{minipage}{0.3 \textwidth}
	\begin{tabular}{r|rrr}
		$\to$: & 0 & 1 & 2 \\
		\hline
		0 & 1 & 1 & 1  \\
		1 & 0 & 1 & 2  \\
		2 & 2 & 1 & 2
	\end{tabular} \hspace{.5cm}	
	
\end{minipage}
\begin{minipage}{0.3 \textwidth}
	\begin{tabular}{r|rrrr}
		$\to$: & 0 & 1 & 2 & 3\\
		\hline
		0 & 1 & 1 & 1 & 1 \\
		1 & 0 & 1 & 2 & 3 \\
		2 & 2 & 1 & 2 & 1 \\
		3 & 3 & 1 & 1 & 3
	\end{tabular}
\end{minipage}
\ \\ \  \\ \ \\
and these are the only simple $\mathbf{DM}$-algebras with more than $2$ elements \cite{balbesDistributive1974}.

\medskip

Hereinafter, we will focus on proving that if $\mathbf A \in \mathbf I$ is a simple algebra with $|A| \geq 3$ then $\mathbf A \in \mathbf{I_{2,0}}$.

\begin{definition} \label{definition_relat_R3}
	Let $\mathbf A \in \mathbf I$. We define the relation $R''$ on $\mathbf A$ as follows:
	$$x \  R^{''} y \mbox{ if and only if } x''=y''.$$
\end{definition}

The following lemma is useful in proving that the relation $R^{''}$ is a congruence on $\mathbf A$.

\begin{Lemma} {\rm \cite[Lemma 3.4]{CoSa2015aI}} \label{Lemma_300315_01}
	Let $\mathbf A$ be an $\mathbf{I}$-zroupoid.  Then $\mathbf A$ satisfies:
	\begin{enumerate} 
		\item[ ] $(x \to y)' \approx (x'' \to y)'$. \label{300315_02}
	\end{enumerate}
\end{Lemma}

\begin{Lemma} \label{lemma_equivalence_relation_R3_simpleAlg}
	Let $\mathbf A \in \mathbf I$. Then $R^{''}$ is a congruence on $\mathbf A$.
\end{Lemma}

\begin{Proof}
	Clearly $R^{''}$ is an equivalence relation. 
	Let $x,y,z,t \in A$ with $x \  R^{''} \ y$ and $z \ R^{''} \  t$. Therefore, $x'' = y''$ and $z'' = t''$. Hence
	by \cite[Lemma 7.5 (b)]{sankappanavarMorgan2012} and Lemma \ref{Lemma_300315_01} 
	we have that $(x \to z)'' = (y \to t)''$.
	Consequently, $R^{''}$ is a congruence on $\mathbf A$.	
\end{Proof}

\begin{Lemma} \label{260315_06}
	Let $\mathbf A \in \mathbf I$. If $\mathbf A \models x'' \approx y''$, then
	\begin{enumerate}[{\rm(a)}]
		\item $\mathbf A \models x'' \approx 0$ \label{260315_01}
		\item $\mathbf A \models (x \to y)' \approx 0$ \label{260315_02}
		\item $\mathbf A \models (x \to y) \to z \approx 0$
	\end{enumerate}
\end{Lemma}

\begin{Proof}
	\begin{enumerate}[(a)]
		\item Let $a \in A$. Then $0 = 0'' = a''$.
		\item Let $a,b \in A$. By (\ref{260315_01}) we have that $a''=b''=0$. Then, from \cite[Lemma 7.5 (b)]{sankappanavarMorgan2012}, $(a \to b)' = (a \to b'')' = (a \to 0)' = a'' = 0$.
		\item Let $a,b,c \in A$. Since (\ref{260315_02}), $0 = [(c' \to a) \to (b \to c)']' = (a \to b) \to c$.
	\end{enumerate}
\end{Proof}

\medskip

The proof of the following lemma is immediate.

\begin{Lemma} \label{260315_10}
	Let $\mathbf A \in \mathbf I$ satisfying the identity $(x \to y)' \approx 0$ then the relation $R^{'}$ defined by,
	$$x R^{'} y \mbox{ if and only if } x' = y',$$
	is a congruence on $\mathbf A$.
\end{Lemma}

\begin{theorem} \label{simpleAlgebras_theorem2}
	Let $\mathbf A \in \mathbf I$ be a simple algebra with $|A| \geq 3$. Then $\mathbf A \in \mathbf{I_{2,0}}$.
\end{theorem}

\begin{Proof}
	Assume that $\mathbf A \not\models x'' \approx x$. Then there exist an element $a \in A$ with
	\begin{equation} \label{260315_03}
	a \not= a''.
	\end{equation}
	Since $0 = 0''$ we conclude that
	\begin{equation} \label{260315_04}
	a \not= 0.
	\end{equation}
	 By Lemma \ref{lemma_equivalence_relation_R3_simpleAlg} we can affirm that the relation 
	 $R^{''}$ (see Definition \ref{definition_relat_R3}) is a congruence on $\mathbf A$.
	By \cite[Corollary 7.7]{sankappanavarMorgan2012}, we have $a'''' = a''$.  Therefore, $a'' R^{''} a$.  Since   $\mathbf A$ is simple, and $a \not= a''$ by \ref{260315_03}, we have
	\begin{equation} \label{260315_05}
	R^{''} = A \times A.
	\end{equation}
	Hence $\mathbf A \models x'' \approx y''$, which,
	by Lemma \ref{260315_06}, implies
	\begin{equation} \label{260315_07}
	\mathbf A \models x'' \approx 0,
	\end{equation}
	\begin{equation} \label{260315_08}
	\mathbf A \models (x \to y)' \approx 0
	\end{equation}
	and
	\begin{equation} \label{260315_09}
	\mathbf A \models (x \to y) \to z \approx 0.
	\end{equation}
	Observe from (\ref{260315_07}) that $a'' = 0 = 0''$, implying $a' R^{'} 0'$.  Also, it follows from (\ref{260315_08}) and Lemma \ref{260315_10} that $R^{'}$ is a congruence.  
	
	If $a' \not= 0'$, since $\mathbf A$ is simple, we would have that $R^{'} = A \times A$. Then $0\ R^{'} a$, and, consequently, $0' = a'$, which is a contradiction. Hence
	\begin{equation} \label{260315_11}
	a' = 0'.
	\end{equation}
	By (\ref{260315_11}), $a R^{'} 0$ and, by (\ref{260315_04}), $R^{'} = A \times A$.
	
	Let $b \in A$. Then $b R^{'} 0'$. We can verify that $b' = 0'' = 0$. Consequently,
	\begin{equation} \label{260315_12}
	\mathbf A \models x' \approx 0.
	\end{equation}
	From (\ref{260315_12}) we conclude that
	$\mathbf A \models 0' \approx 0.$ Then, for $b \in A$, we have $0 \to b = 0' \to b = (0 \to 0) \to b = 0$ by (\ref{260315_09}).
	Hence
	\begin{equation} \label{260315_13}
	\mathbf A \models 0 \to x \approx 0.
	\end{equation}
	Now we consider the following cases:
	\begin{itemize}
		\item Assume that $\mathbf A \models x \to y \approx 0$. Let $R_a$ denote the equivalence relation correspnding to the partition $\{\{0\}, A \setminus \{0\}\}$ of $A$. Notice that $|A \setminus \{0\}| \geq 2$ since $|A| \geq 3$, whence $R^{'}$ is nontrivial.  Since $\mathbf A \models x \to y \approx 0$, $R_a$ is a congruence, which is a contradiction since the algebra is simple.
		\item Assume that $\mathbf A \not\models x \to y \approx 0$. Then the set $$T = \{b \in A:\ \mbox{ exists } c \in A \mbox{ with } b \to c \not= 0 \}$$ is not empty, and from (\ref{260315_13}), $0 \not\in T$. Then the equivalence relation $R_T$ associated to the partition $\{T, A \setminus T\}$ of $A$ satisfies that $R_T \not= \Delta$ and $R_T \not= A \times A$. Let the elements $b,c,d,e, f$ in $A$ with $b \ R_T \ c$ and $d \ R_T \ e$.  Then, from (\ref{260315_09}), we have that $(b \to d) \to f = (c \to e) \to f = 0$. Hence $b \to d \ R_T \ c \to e$, implying that $R_T$ is a congruence on $\mathbf A$.  Since $R_T \not= \Delta$, we have a contradiction, because $\mathbf A$ is simple and $|A| \geq 3$.
	\end{itemize}
	
Thus we have proved that $\mathbf A \models x'' \approx x$, completing the proof.
\end{Proof}

\ \\ 

In view of Theorem \ref{simpleAlgebras_theorem1} and Theorem \ref{simpleAlgebras_theorem2} we have the following crucial result.

\begin{corollary}
	Let $\mathbf A \in \mathbf I$ be a simple algebra with $|A| \geq 3$. Then $\mathbf A \in \mathbf{DM}$.
\end{corollary}

Thus we have proved the following desired result describing the simple algebras of $\mathbf I$.

\begin{theorem}
	The only simple algebras in $\mathbf I$ are $\mathbf{2_z}$, $\mathbf{2_s}$, $\mathbf{2_b}$, $\mathbf 3_k$ and $\mathbf 4_d$.
\end{theorem}

Recall (\cite{burrisCourse1981}) that a variety is semisimple if and only if every subdirectly irreducible algebra in it is simple. 
\begin{corollary}
A subvariety $\mathbf{V}$ of $\mathbf{I}$ is semi-simple if and only if $\mathbf V \subseteq \mathbb V(\mathbf{2_z}, \mathbf{2_s},  \mathbf 4_d)$.	
\end{corollary}

\begin{Corollary} \label{latticeSemiSimpleVarieties}
	The lattice of semisimple subvarieties of $\mathbf I$ is isomorphic to the direct product of a 4-element Boolean lattice and and a 4-element chain and is  
	the one depicted in the following figure, where $\mathbf A_1, \ldots, \mathbf A_n$ denotes $\mathbb V(\mathbf A_1, \ldots, \mathbf A_n)$.
\end{Corollary}

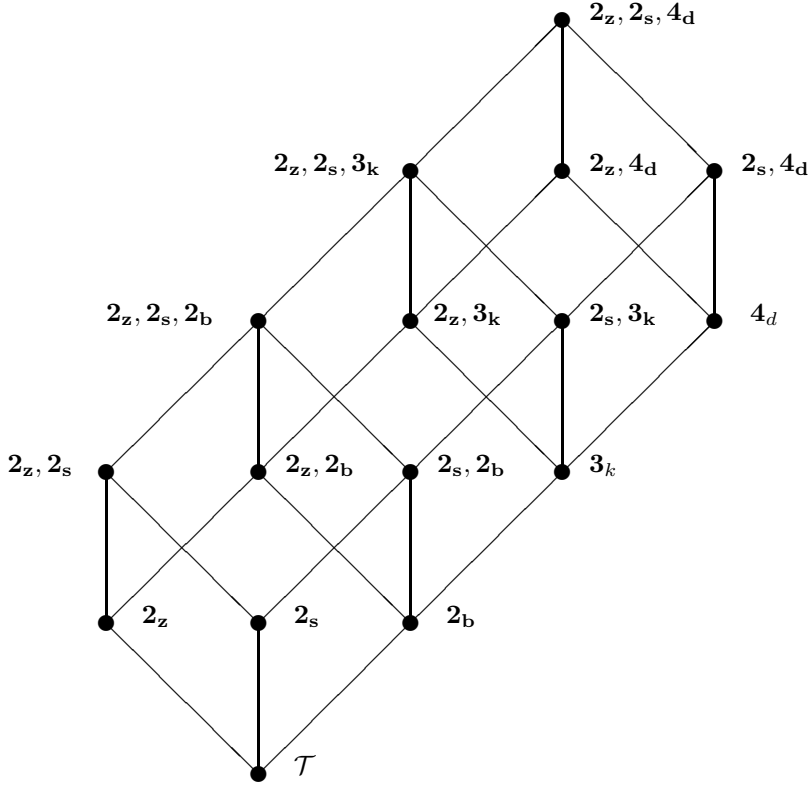
\begin{small}

	\begin{center}
		
		\setlength{\unitlength}{1cm}
		\begin{picture}(14,14)(0,0)  \label{graficoISSH}
		\linethickness{0.02cm}
		
		
		\put(4,2){\circle*{0.2}}
		
		\put(2,4){\circle*{0.2}}
		\put(4,4){\circle*{0.2}}
		\put(6,4){\circle*{0.2}}
		
		\put(2,6){\circle*{0.2}}
		\put(4,6){\circle*{0.2}}
		\put(6,6){\circle*{0.2}}
		\put(8,6){\circle*{0.2}}
		
		\put(4,8){\circle*{0.2}}
		\put(6,8){\circle*{0.2}}
		\put(8,8){\circle*{0.2}}
		\put(10,8){\circle*{0.2}}
		
		\put(6,10){\circle*{0.2}}
		\put(8,10){\circle*{0.2}}
		\put(10,10){\circle*{0.2}}
		
		\put(8,12){\circle*{0.2}}
		
		
		\put(4,2){$\ \ \ \ \mathcal T$}
		
		\put(2,4){$\ \ \ \  \mathbf{2_z}$}
		\put(4,4){$\ \ \ \  \mathbf{2_s}$}
		\put(6,4){$\ \ \ \  \mathbf{2_b}$}
		
		\put(0.7,6){$\mathbf{2_z}, \mathbf{2_s}$}
		\put(4,6){$\ \ \  \mathbf{2_z}, \mathbf{2_b}$}
		\put(6,6){$\ \ \  \mathbf{2_s}, \mathbf{2_b}$}
		\put(8,6){$\ \ \  \mathbf 3_k$}
		
		\put(2,8){$\mathbf{2_z}, \mathbf{2_s}, \mathbf{2_b}$}
		\put(6.3,8){$\mathbf{2_z},  \mathbf{3_k}$}
		\put(8,8){$\ \ \  \mathbf{2_s}, \mathbf{3_k}$}
		\put(10,8){$\ \ \ \  \mathbf 4_d$}
		
		\put(4.2,10){$\mathbf{2_z}, \mathbf{2_s}, \mathbf{3_k}$}
		\put(8,10){$\ \ \  \mathbf{2_z},  \mathbf{4_d}$}
		\put(10,10){$\ \ \  \mathbf{2_s},  \mathbf{4_d}$}
		
		\put(8,12){$\ \ \ \mathbf{2_z}, \mathbf{2_s}, \mathbf{4_d}$}
		
		
		
		\put(4,2){\line(0,1){2}}
		
		\put(2,4){\line(0,1){2}}
		\put(6,4){\line(0,1){2}}
		
		\put(4,6){\line(0,1){2}}
		\put(8,6){\line(0,1){2}}
		
		\put(6,8){\line(0,1){2}}
		\put(10,8){\line(0,1){2}}
		
		\put(8,10){\line(0,1){2}}

		
		\put(4,2){\line(1,1){6}}
		\put(2,4){\line(1,1){6}}
		\put(4,4){\line(1,1){6}}
		\put(2,6){\line(1,1){6}}
		
		\put(2,4){\line(1,-1){2}}
		\put(2,6){\line(1,-1){2}}
		
		\put(4,6){\line(1,-1){2}}
		\put(4,8){\line(1,-1){2}}
		
		\put(6,8){\line(1,-1){2}}
		\put(6,10){\line(1,-1){2}}
		
		\put(8,10){\line(1,-1){2}}
		\put(8,12){\line(1,-1){2}}
		
		\end{picture}
		
	\end{center}
	
\end{small}

\begin{Corollary}
Semisimple varieties of {\bf{I}}-zroupoids are locally finite.
\end{Corollary}

We conclude this section by raising the problem of  axiomatization of subvarieties of  $\mathbb V(\mathbf{2_z}, \mathbf{2_s},  \mathbf 4_d)$.

\section{Appendix}

\begin{Proof} {\bf of Lemma \ref{general_properties2}}
Items from (\ref{281014_05}) to (\ref{071114_02}) are proved in \cite{CoSa2015aI}.
	Let $a,b,c,d \in A$.
	\begin{enumerate}
		\item[(\ref{291014_08})]
		$$
		\begin{array}{lcll}
		[(0 \to a) \to b] \to a  & = & [(b \to a) \to (0 \to a)']' & \mbox{by (\ref{291014_03})} \\
		& = & (b \to a)'' & \mbox{by (\ref{291014_07}) using } x = b, y = a \\
		& = & b \to a
		\end{array}
		$$
		
		\item[(\ref{311014_06})]
		$$
		\begin{array}{lcll}
		a \to (0 \to b) & = & [\{0 \to (0 \to b)\} \to a] \to (0 \to b)  & \mbox{by (\ref{291014_08}) with } x = 0 \to b, y = a  \\
		& = & [(0 \to b) \to a] \to (0 \to b) & \mbox{by (\ref{311014_03})} \\
		& = & (0 \to a) \to (0 \to b) & \mbox{by (\ref{291014_10})}
		\end{array}
		$$
		
		\item[(\ref{031114_03})]
		$$
		\begin{array}{lcll}
		(a \to b)' \to (0 \to c) & = & [(a \to b) \to 0] \to (0 \to c) & \mbox{} \\
		& = & [\{(0 \to c)' \to (a \to b)\} \to \{0 \to (0 \to c)\}']' & \mbox{by (I)} \\
		& = & [\{(0 \to c)' \to (a \to b)\} \to (0 \to c)']' & \mbox{by (\ref{311014_03})} \\
		& = & [\{0 \to (a \to b)\} \to (0 \to c)']' & \mbox{by (\ref{291014_10}) taking } x = a \to b,\\
		&    &   & y = (0 \to c)' \\
		& = & [\{0 \to (a \to b)\} \to (c' \to 0')']' & \mbox{by Lemma \ref{general_properties} (\ref{cuasiConmutativeOfImplic2})} \\
		& = & [(a \to b) \to c'] \to 0' & \mbox{from (I)} \\
		& = & 0 \to [(a \to b) \to c']' & \mbox{by Lemma \ref{general_properties} (\ref{cuasiConmutativeOfImplic2})}
		\end{array}
		$$
		
		\item[(\ref{311014_01})]  
		$$
		\begin{array}{lcll}
		[(a' \to b') \to (b \to a)'] \to (b \to a) & = & [(a' \to b') \to (b \to a)']'' \to (b \to a) & \mbox{} \\
		& = & [(b' \to b) \to a]' \to (b \to a) & \mbox{by (I)} \\
		& = & [b \to a]' \to (b \to a) & \mbox{by Lemma \ref{general_properties_equiv} (\ref{LeftImplicationwithtilde})} \\
		& = & b \to a & \mbox{by Lemma \ref{general_properties_equiv} (\ref{LeftImplicationwithtilde})}
		\end{array}
		$$
		
		\item[(\ref{311014_02})] 
		$$
		\begin{array}{lcll}
		[0 \to (a' \to b')'] \to (b \to a)  & = & [(a' \to b') \to 0'] \to (b \to a) & \mbox{by Lemma \ref{general_properties} (\ref{cuasiConmutativeOfImplic2})} \\
		& = & [(a' \to b') \to (b \to a)'] \to (b \to a) & \mbox{by (\ref{281014_05})} \\
		& = & b \to a & \mbox{from (\ref{311014_01})}
		\end{array}
		$$
		
		\item[(\ref{311014_04})]
		$$
		\begin{array}{lcll}
		(0 \to (a \to b)) \to a & = & [(a' \to 0) \to \{(a \to b) \to a\}']' & \mbox{from (I)} \\
		& = & [(a' \to 0) \to \{(0 \to b) \to a\}']' & \mbox{by (\ref{291014_10})} \\
		& = & [\{0 \to (0 \to b\}] \to a & \mbox{using (I)} \\
		& = & (0 \to b) \to a & \mbox{by (\ref{311014_03})}
		\end{array}
		$$
		
		\item[(\ref{311014_05})]   
		$$
		\begin{array}{lcll}
		(0 \to a) \to (a \to b)  & = & [(a \to b) \to a] \to (a \to b) & \mbox{by (\ref{291014_10}) with } y = a \to b, x = a \\
		& = & [a \to (b \to a)']' \to (a \to b) & \mbox{by (\ref{291014_09})} \\
		& = & [a'' \to (b \to a)']' \to (a \to b) & \mbox{} \\
		& = & [(a' \to 0) \to (b \to a)']' \to (a \to b) & \mbox{} \\
		& = & [(0 \to b) \to a]  \to (a \to b) & \mbox{from (I)} \\
		& = & [\{0 \to (a \to b)\} \to a]  \to (a \to b) & \mbox{by (\ref{311014_04})} \\
		& = & a \to (a \to b) & \mbox{by (\ref{291014_08}) with } x = a \to b, y = a
		\end{array}
		$$
		
		\item[(\ref{311014_07})]   
		$$
		\begin{array}{lcll}
		[b \to (0 \to a)]' & = & [(0 \to b) \to (0 \to a)]' & \mbox{by (\ref{311014_06})} \\
		& = & [\{(0 \to a)' \to 0\} \to \{b \to (0 \to a)\}']'' & \mbox{from (I)} \\
		& = & [(0 \to a)' \to 0] \to [b \to (0 \to a)]' & \mbox{} \\
		& = & (0 \to a)'' \to [b \to (0 \to a)]' & \mbox{} \\
		& = & (0 \to a) \to [b \to (0 \to a)]' & \mbox{}
		\end{array}
		$$
		
		\item[(\ref{311014_08})]
		$$
		\begin{array}{lcll}
		[(a \to b) \to (0 \to a)]' & = & [\{(0 \to a)' \to a\} \to \{b \to (0 \to a)\}']'' & \mbox{from (I)} \\
		& = &[(0 \to a)' \to a] \to [b \to (0 \to a)]'  & \mbox{} \\
		& = & [\{(0 \to a) \to 0\} \to a] \to [b \to (0 \to a)]' & \mbox{} \\
		& = & (0 \to a) \to [b \to (0 \to a)]' & \mbox{by (\ref{291014_08}) with } x = a, y = 0
		\end{array}
		$$
		
		\item[(\ref{031114_01})]  
		$$
		\begin{array}{lcll}
		b \to (0 \to a) & = & [b \to (0 \to a)]'' & \mbox{} \\
		& = & [(0 \to a) \to \{b \to (0 \to a)\}']' & \mbox{by (\ref{311014_07})} \\
		& = & [(a \to b) \to (0 \to a)]'' & \mbox{by (\ref{311014_08})} \\
		& = & (a \to b) \to (0 \to a). & \mbox{}
		\end{array}
		$$
		
		\item[(\ref{031114_02})]
		$$
		\begin{array}{lcll}
		b \to (b \to a)  & = & (0 \to b) \to (b \to a) & \mbox{by (\ref{311014_05})} \\
		& = & [0 \to (0 \to b)] \to (b \to a) & \mbox{by (\ref{311014_03})} \\
		& = & [(b \to a) \to (0 \to b)] \to (b \to a) & \mbox{by (\ref{291014_10})} \\
		& = & [a \to (0 \to b)] \to (b \to a) & \mbox{by (\ref{031114_01})} \end{array}
		$$
		
		\item[(\ref{031114_04})]   
		$$
		\begin{array}{lcll}
		a \to b & = & [0 \to (b' \to a')'] \to (a \to b) & \mbox{by (\ref{311014_02})} \\
		& = & [0 \to \{(b \to 0) \to a'\}'] \to (a \to b) & \mbox{} \\
		& = & [(b \to 0)' \to (0 \to a)] \to (a \to b) & \mbox{by (\ref{031114_03})} \\
		& = & [b'' \to (0 \to a)] \to (a \to b) & \mbox{} \\
		& = & [b \to (0 \to a)] \to (a \to b) & \mbox{} \\
		& = & a \to (a \to b) & \mbox{by (\ref{031114_02})}
		\end{array}
		$$
		
		\item[(\ref{171114_01})]
		$$
		\begin{array}{lcll}
		c \to ((a \to b) \to c)' & = & [\{c \to (a \to b)\} \to c]' & \mbox{by (\ref{291014_09})} \\
		& = & [\{0 \to (a \to b)\} \to c]' & \mbox{by (\ref{291014_10})} \\
		& = & [\{a \to (0 \to b)\} \to c]' & \mbox{by (\ref{071114_04})}
		\end{array}
		$$
		
		\item[(\ref{181114_01})]  
		$$
		\begin{array}{lcll}
		[0 \to (a \to b)'] \to b & = & [(a \to b) \to 0'] \to b & \mbox{by Lemma \ref{general_properties} (\ref{cuasiConmutativeOfImplic2})} \\
		& = & [(a \to b) \to b'] \to b  & \mbox{by (\ref{281014_05})} \\
		& = & [(a \to 0') \to b'] \to b & \mbox{by (\ref{281014_05})} \\
		& = & [(a \to 0') \to 0'] \to b & \mbox{by (\ref{281014_05})} \\
		& = & [(0 \to a') \to 0'] \to b  & \mbox{by Lemma \ref{general_properties} (\ref{cuasiConmutativeOfImplic2})} \\
		& = & [0 \to (0 \to a')'] \to b  & \mbox{by Lemma \ref{general_properties} (\ref{cuasiConmutativeOfImplic2})} \\
		& = & (0 \to a'') \to b & \mbox{by (\ref{031114_07})} \\
		& = & (0 \to a) \to b & \mbox{}
		\end{array}
		$$
		
		\item[(\ref{181114_02})]   
		$$
		\begin{array}{lcll}
		[\{a \to (0 \to (b \to c)')\} \to c]' & = & (c' \to a) \to [\{0 \to (b \to c)'\} \to c]' & \mbox{by (I)} \\
		& = & (c' \to a) \to [(0 \to b) \to c]' & \mbox{by (\ref{181114_01})} \\
		& = & (c' \to a) \to [(c \to b) \to c]' & \mbox{by (\ref{291014_10})} \\
		& = & (c' \to a) \to [c \to (b \to c)']'' & \mbox{by (\ref{291014_09})} \\
		& = & (c' \to a) \to [c \to (b \to c)'] & \mbox{}
		\end{array}
		$$
		
		\item[(\ref{181114_03})]
		$$
		\begin{array}{lcll}
		a \to [\{b \to (c \to a)'\} \to a]' & = & [\{b \to (0 \to (c \to a)')\} \to a]' & \mbox{by (\ref{171114_01}) with } x = b, y = (c \to a)',\\
		&  &  & z = a \\
		& = & (a' \to b) \to \{a \to (c \to a)'\} & \mbox{by (\ref{181114_02}) with } x = b, y = c, z = a
		\end{array}
		$$
		
		\item[(\ref{181114_11})]  
		$$
		\begin{array}{lcll}
		\{0 \to (a \to b)\} \to b' & = & [(b \to 0) \to \{(a \to b) \to b'\}']' & \mbox{from (I)} \\
		& = & [(b \to 0) \to \{(a \to 0') \to b'\}']' & \mbox{by (\ref{281014_05})} \\
		& = & \{0 \to (a \to 0')\} \to b' &  \mbox{from (I)} \\
		& = & \{0 \to (0 \to a')\} \to b'  & \mbox{by Lemma \ref{general_properties} (\ref{cuasiConmutativeOfImplic2})} \\
		& = & (0 \to a') \to b' & \mbox{by (\ref{311014_03})} \\
		& = & (a \to 0') \to b' & \mbox{by Lemma \ref{general_properties} (\ref{cuasiConmutativeOfImplic2})} \\
		& = & (a \to b'') \to b' & \mbox{by (\ref{281014_05})} \\
		& = & (a \to b) \to b' & \mbox{} \\
		& = & b \to (a \to b)' & \mbox{by (\ref{071114_05})}
		\end{array}
		$$
		
		\item[(\ref{181114_12})]  
		$$
		\begin{array}{lcll}
		(a' \to b) \to \{a \to (c \to a)'\} & = & a \to [\{b \to (c \to a)'\} \to a]' & \mbox{by (\ref{181114_03})} \\
		& = & (0 \to b) \to [\{0 \to (c \to a)\} \to a'] & \mbox{by (\ref{181114_10}) with } x = a, y = b, \\
		&    &   & z = c \to a \\
		& = & (0 \to b) \to [a \to (c \to a)'] & \mbox{by (\ref{181114_11}) with } x = c, y = a
		\end{array}
		$$
		
		\item[(\ref{181114_13})]  
		$$
		\begin{array}{lcll}
		a \to [(c \to b) \to a]'  & = & [\{c \to (0 \to b)\} \to a]' & \mbox{by (\ref{171114_01}) with } x = c, y = b, z = a \\
		& = & [\{(b \to c) \to (0 \to b)\} \to a]' & \mbox{by (\ref{031114_01})} \\
		& = & (a' \to (b \to c)) \to [(0 \to b) \to a]' & \mbox{by (I)} \\
		& = & \{a' \to (b \to c)\} \to [(a \to b) \to a]' & \mbox{by (\ref{291014_10})} \\
		& = & \{a' \to (b \to c)\} \to [a \to (b \to a)']'' & \mbox{by (\ref{291014_09})} \\
		& = & \{a' \to (b \to c)\} \to [a \to (b \to a)'] & \mbox{}
		\end{array}
		$$
		
		\item[(\ref{181114_14})]
		$$
		\begin{array}{lcll}
		a \to [(b \to c) \to a]' & = & [a' \to (c \to b)] \to [a \to (c \to a)'] & \mbox{by (\ref{181114_13}) with } x = a, y = c, z = b \\
		& = & [0 \to (c \to b)] \to [a \to (c \to a)'] & \mbox{by (\ref{181114_12}) with } x = a, y = c \to b, z = c
		\end{array}
		$$
		
		\item[(\ref{181114_15})]  
		$$
		\begin{array}{lcll}
		[0 \to (a \to b)] \to [c \to (a \to c)'] & = & c \to [(b \to a) \to c]' & \mbox{by (\ref{181114_14}) with } x = c, y = b, z = a \\
		& = & [\{c \to (b \to a)\} \to c]' & \mbox{by (\ref{291014_09})} \\
		& = & [\{0 \to (b \to a)\} \to c]' & \mbox{by (\ref{291014_10})} \\
		& = & [\{b \to (0 \to a)\} \to c]' & \mbox{by (\ref{071114_04})} \\
		& = & (c' \to b) \to [(0 \to a) \to c]' & \mbox{from (I)} \\
		& = & (c' \to b) \to [(c \to a) \to c]' & \mbox{by (\ref{291014_10})} \\
		& = & (c' \to b) \to [c \to (a \to c)'] & \mbox{by (\ref{291014_09})} \\
		& = & (0 \to b) \to [c \to (a \to c)'] & \mbox{by (\ref{181114_12}) with } x = c, y = b, z = a
		\end{array}
		$$
		
		\item[(\ref{181114_16})]
		$$
		\begin{array}{lcll}
		a \to [(b \to c) \to a]'  & = & [0 \to (c \to b)] \to [a \to (c \to a)'] & \mbox{by (\ref{181114_14})} \\
		& = & (0 \to b) \to [a \to (c \to a)'] & \mbox{by (\ref{181114_15}) with } x = c, y = b, z = a
		\end{array}
		$$
		
		\item[(\ref{191114_02})] This follows immediately from (\ref{071114_03}).
		
		\item[(\ref{031214_16})]
		$$
		\begin{array}{lcll}
		a \to (b \to a)' & = & [(a \to b) \to a]' & \mbox{by (\ref{291014_09})} \\
		& = & [(0 \to b) \to a]' & \mbox{by (\ref{291014_10})} \\
		& = & [(b' \to 0') \to a]' & \mbox{by Lemma \ref{general_properties} (\ref{cuasiConmutativeOfImplic2})} \\
		& = & (0 \to b') \to a' & \mbox{by (\ref{071114_02})} \\
		& = & (b \to 0') \to a' & \mbox{by Lemma \ref{general_properties} (\ref{cuasiConmutativeOfImplic2})}
		\end{array}
		$$

		\item[(\ref{130315_02})]
		$$
		\begin{array}{lcll}
		(a' \to b) \to a'  & = & [a' \to (b \to a')']' & \mbox{by  (\ref{291014_09})} \\
		& = & [(b \to 0') \to a'']' & \mbox{by  (\ref{031214_16})} \\
		& = & [(b \to 0') \to a]' & \mbox{} \\
		& = & [(0 \to b') \to a]' & \mbox{by Lemma \ref{general_properties} (\ref{cuasiConmutativeOfImplic})} \\
		& = & [(a \to b') \to a]' & \mbox{by  (\ref{291014_10})} \\
		& = & a' & \mbox{by Hyphotesis}
		\end{array}
		$$
		
		\item[(\ref{130315_08})]   
		$$
		\begin{array}{lcll}
		[\{b \to (0 \to c)\} \to (a \to 0')']' & = & [\{(a \to 0')'' \to b\} \to \{(0 \to c) \to (a \to 0')'\}']'' & \mbox{by (I)} \\
		& = & [(a \to 0') \to b] \to [(0 \to c) \to (a \to 0')']' & \mbox{} \\
		& = & [(a \to 0') \to b] \to [(c \to a) \to 0'] & \mbox{by (I)} \\
		& = & [(a \to 0') \to b] \to [0 \to (c \to a)'] & \mbox{by Lemma \ref{general_properties} (\ref{cuasiConmutativeOfImplic})} \\
		& = & [(0 \to a') \to b] \to [0 \to (c \to a)'] & \mbox{by Lemma \ref{general_properties} (\ref{cuasiConmutativeOfImplic})}
		\end{array}
		$$
		
		\item[(\ref{130315_09})]   
		$$
		\begin{array}{lcll}
		[\{a \to (0 \to b)\} \to c]'& = & [(c' \to a) \to \{(0 \to b) \to c\}']'' & \mbox{by (I)} \\
		& = & (c' \to a) \to [(0 \to b) \to c]' & \mbox{} \\
		& = & (c' \to a) \to [(c \to b) \to c]' & \mbox{by  (\ref{291014_10})} \\
		& = & (c' \to a) \to [c \to (b \to c)']'' & \mbox{by  (\ref{291014_09})} \\
		& = & (c' \to a) \to [c \to (b \to c)'] & \mbox{} \\
		& = & (0 \to a) \to [c \to (b \to c)'] & \mbox{by  (\ref{181114_12})}
		\end{array}
		$$

		\item[(\ref{130315_10})]

		$$
		\begin{array}{lcll}
		[(0 \to a') \to b] \to [0 \to (c \to a)'] & = & [\{b \to (0 \to c)\} \to (a \to 0')']' & \mbox{by (\ref{130315_08})}  \\
		& = & (0 \to b) \to [(a \to 0')' \to \{c \to (a \to 0')'\}'] & \mbox{by (\ref{130315_09}) with }  \\
		&  &  & x = b, y = c, \\
		&   &   & z = (a \to 0')'   \\
		& = & (0 \to b) \to [(a \to 0')' \to \{c \to (a \to 0')'\}']'' & \mbox{} \\
		& = & (0 \to b) \to [\{(a \to 0')' \to c\} \to (a \to 0')']' & \mbox{by (\ref{291014_09})} \\
		& = & (0 \to b) \to [(0 \to c) \to (a \to 0')']' & \mbox{by (\ref{291014_10})} \\
		& = &(0 \to b) \to [(c \to a) \to 0']  & \mbox{by (I)} \\
		& = & (0 \to b) \to [0 \to (c \to a)'] & \mbox{by Lemma \ref{general_properties} (\ref{cuasiConmutativeOfImplic})} \\
		& = & b \to [0 \to (c \to a)'] & \mbox{by (\ref{311014_06})}
		\end{array}
		$$
		
		\item[(\ref{130315_11})]  
		
		$$
		\begin{array}{lcll}
		[(0 \to a) \to b] \to (0 \to c) & = & [0 \to \{(0 \to a) \to b\}] \to (0 \to c) & \mbox{by (\ref{311014_06})} \\
		& = & [(0 \to a) \to (0 \to b)] \to (0 \to c) & \mbox{by  (\ref{071114_04})} \\
		& = & [a \to (0 \to b)] \to (0 \to c) & \mbox{by (\ref{311014_06})} \\
		& = & [0 \to (a \to b)] \to (0 \to c) & \mbox{by (\ref{071114_04})} \\
		& = & (a \to b) \to (0 \to c) & \mbox{by (\ref{311014_06})} \\
		& = & 0 \to ((a \to b) \to c) & \mbox{by (\ref{071114_04})}
		\end{array}
		$$
		
		\item[(\ref{130315_12})]
		$$
		\begin{array}{lcll}
		b \to [0 \to (c \to a)'] & = & [(0 \to a') \to b] \to [0 \to (c \to a)']  & \mbox{by (\ref{130315_10})} \\
		& = & 0 \to [(a' \to b) \to (c \to a)'] & \mbox{by (\ref{130315_11})} \\
		& = & (a' \to b) \to [0 \to (c \to a)'] & \mbox{by (\ref{071114_04})}
		\end{array}
		$$
		
		\item[(\ref{130315_13})]
		$$
		\begin{array}{lcll}
		b \to [0 \to (c \to a)'] & = & (a' \to b) \to [0 \to (c \to a)'] & \mbox{by (\ref{130315_12})} \\
		& = & 0 \to [(a' \to b) \to (c \to a)'] & \mbox{by (\ref{071114_04})}
		\end{array}
		$$

		\item[(\ref{140315_01})]  
		$$
		\begin{array}{lcll}
		(0 \to a) \to  [\{0 \to (b \to c)\} \to d'] & = & (0 \to a) \to  [\{0 \to (b \to c)\} \to \{d \to (0 \to d)'\}] & \mbox{by (\ref{291014_02})} \\
		& = & (0 \to a) \to [d \to [(b \to c)' \to d]]' & \mbox{by (\ref{181114_16}) with }  \\
		&  &  & x = d, y = b \to c, \\
		&   &  &  z = 0   \\
		& = & [\{a \to (0 \to (b \to c)')\}\to d]' & \mbox{by (\ref{130315_09}) with }  \\
		&  &  & x = a, y = (b \to c)',\\
		&   &  & z = d   \\
		& = &[[0 \to \{(c' \to a) \to (b \to c)'\}] \to d]'  & \mbox{by (\ref{130315_13}) with }  \\
		&  &  & x = c, y = a, z = b   \\
		& = & [[0 \to \{(c' \to a) \to (b \to c)'\}''] \to d]' & \mbox{} \\
		& = & [[0 \to \{(a \to b) \to c\}'] \to d]' & \mbox{by (I)} \\
		& = & [[d \to [(a \to b) \to c]'] \to d]' & \mbox{by (\ref{291014_10})} \\
		& = & [d \to \{((a \to b) \to c)' \to d\}']'' & \mbox{by (\ref{291014_09})} \\
		& = & d \to [\{(a \to b) \to c\}' \to d]' & \mbox{} \\
		& = & [0 \to \{(a \to b) \to c\}] \to [d \to (0 \to d)'] & \mbox{by (\ref{181114_16}) with } x = d,  \\
		&  &   &  y = (a \to b) \to c,\\
		&   &   &  z = 0 \\
		& = & [0 \to \{(a \to b) \to c\}] \to d' & \mbox{by (\ref{291014_02})}
		\end{array}
		$$
		
		\item[(\ref{140315_02})]

		$$
		\begin{array}{lcll}
		(0 \to a) \to [(0 \to b') \to (0 \to c)'] & = & (0 \to a) \to [\{0 \to (b \to 0)\} \to (0 \to c)'] & \mbox{} \\
		& = & [0 \to \{(a \to b) \to 0\}] \to (0 \to c)' & \mbox{by (\ref{140315_01}) with }   \\
		&  &  & x = a, y = b, z = 0,\\
		&   &  &  u = 0 \to c   \\
		& = & [0 \to (a \to b)'] \to (0 \to c)' & \mbox{} \\
		& = & [(a \to b) \to (0 \to c)]' & \mbox{by (\ref{031114_06})} \\
		& = & [0 \to \{(a \to b) \to c\}]' & \mbox{by (\ref{071114_04})}
		\end{array}
		$$
		
		\item[(\ref{140315_03})]
		
		$$
		\begin{array}{lcll}
		0 \to [(a \to b) \to (c \to a')'] & = & 0 \to [(a'' \to b) \to (c \to a')']''  & \mbox{} \\
		& = & 0 \to [(b \to c) \to a']' & \mbox{by (I)} \\
		& = & 0 \to [0 \to \{(b \to c) \to a'\}'] & \mbox{by (\ref{311014_03})} \\
		& = & 0 \to [(0 \to b) \to \{(0 \to c') \to (0 \to a')'\}] & \mbox{by (\ref{140315_02}) with }  \\
		&  &  & x = b, y = c, z = a'   \\
		& = & 0 \to [(0 \to b) \to \{(0 \to c') \to (a \to 0')'\}] & \mbox{by Lemma \ref{general_properties} (\ref{cuasiConmutativeOfImplic})} \\
		& = & 0 \to [(0 \to b) \to \{(0 \to c') \to (a \to 0')'\}''] & \mbox{} \\
		& = & 0 \to [(0 \to b) \to \{(c' \to a) \to 0'\}'] & \mbox{by (I)} \\
		& = & 0 \to [(0 \to b) \to \{(c' \to a) \to 0'\}']'' & \mbox{} \\
		& = & 0 \to [\{b \to (c' \to a)\} \to 0']' & \mbox{by (I)} \\
		& = & 0 \to [0 \to \{b \to (c' \to a)\}']' & \mbox{by Lemma \ref{general_properties} (\ref{cuasiConmutativeOfImplic})} \\
		& = & 0 \to [b \to (c' \to a)]'' & \mbox{by (\ref{031114_07})} \\
		& = & 0 \to [b \to (c' \to a)] & \mbox{}
		\end{array}
		$$
		
		\item[(\ref{140315_04})]
		
		$$
		\begin{array}{lcll}
		0 \to [(a \to b)' \to c] & = & 0 \to [(a \to b) \to c']' & \mbox{by (\ref{191114_05})} \\
		& = & 0 \to [(c'' \to a) \to (b \to c')']'' & \mbox{by (I)} \\
		& = & 0 \to [(c \to a) \to (b \to c')'] & \mbox{} \\
		& = & 0 \to [a \to (b' \to c)] & \mbox{by (\ref{140315_03})}
		\end{array}
		$$

		\item[(\ref{140315_05})]  
		
		$$
		\begin{array}{lcll}
		[\{(a \to b) \to c\} \to \{0 \to (b \to c)\}'] \to \{(a \to b) \to c\} & = & [0 \to \{0 \to (b \to c)\}'] \to \{(a \to b) \to c\} & \mbox{by (\ref{291014_10})} \\
		& = & [0 \to (b \to c)'] \to [(a \to b) \to c] & \mbox{by (\ref{311014_03})} \\
		& = & [0 \to (b \to c)'] \to [(c' \to a) \to (b \to c)']' & \mbox{by (I)} \\
		& = & [(c' \to a) \to (b \to c)']' & \mbox{by (\ref{291014_06})} \\
		& = & (a \to b) \to c  & \mbox{by (I)}
		\end{array}
		$$
		
		\item[(\ref{170315_01})]
		$$
		\begin{array}{lcll}
		[0 \to (a \to b)'] \to (c \to b) & = & [0 \to (a'' \to b)'] \to (c \to b) & \mbox{} \\
		& = & [a' \to (0 \to b')] \to (c \to b) & \mbox{by (\ref{071114_01})} \\
		& = & [0 \to \{(0 \to a') \to b'\}] \to (c \to b) & \mbox{by (\ref{191114_02})} \\
		& = & [(0 \to a') \to (0 \to b')] \to (c \to b) & \mbox{by (\ref{071114_04})} \\
		& = & [\{(c \to b)' \to (0 \to a')\} \to \{(0 \to b') \to (c \to b)\}']' & \mbox{by (I)} \\
		& = & [\{(c \to b)' \to (0 \to a')\} \to (c \to b)']' & \mbox{by (\ref{281014_07})} \\
		& = & [(c \to b)' \to \{(0 \to a') \to (c \to b)'\}']'' & \mbox{by (\ref{291014_09})} \\
		& = & (c \to b)' \to [(0 \to a') \to (c \to b)']' & \mbox{} \\
		& = & [(0 \to a') \to 0'] \to (c \to b)'' & \mbox{by (\ref{031214_16})} \\
		& = & [(0 \to a') \to 0'] \to (c \to b) & \mbox{} \\
		& = & [0 \to (0 \to a')'] \to (c \to b) & \mbox{by Lemma \ref{general_properties} (\ref{cuasiConmutativeOfImplic})} \\
		& = & (0 \to a'') \to (c \to b)  & \mbox{by (\ref{031114_07})} \\
		& = & (0 \to a) \to (c \to b) & \mbox{}
		\end{array}
		$$

		\item[(\ref{250315_01})]
		
		$$
		\begin{array}{lcll}
		[b \to [0 \to (0 \to c)']] \to a & = & [b \to (0 \to c')] \to a & \mbox{by (\ref{031114_07})} \\
		& = & [0 \to (b' \to c)'] \to a & \mbox{by (\ref{071114_01})} \\
		& = & [a \to (b' \to c)'] \to a & \mbox{by (\ref{291014_10})}
		\end{array}
		$$
		
		\item[(\ref{250315_02})]
		$$
		\begin{array}{lcll}
		b \to (c \to a) & = & [\{0 \to (c \to a)\} \to b] \to (c \to a) & \mbox{by (\ref{291014_08})} \\
		& = & [[(c \to a)' \to [0 \to (c \to a)]] \to [b \to (c \to a)]']' & \mbox{by (I)} \\
		& = & [[(c \to a)' \to [(c \to a)' \to 0']] \to [b \to (c \to a)]']' & \mbox{by Lemma \ref{general_properties} (\ref{cuasiConmutativeOfImplic})} \\
		& = & [[(c \to a)' \to 0'] \to [b \to (c \to a)]']' & \mbox{by (\ref{031114_04})} \\
		& = & [[0 \to (c \to a)] \to [b \to (c \to a)]']' & \mbox{by Lemma \ref{general_properties} (\ref{cuasiConmutativeOfImplic})} \\
		& = & [[c \to (0 \to a)] \to [b \to (c \to a)]']' & \mbox{by (\ref{071114_04})} \\
		& = & [[c'' \to (0 \to a)]'' \to [b \to (c \to a)]']' & \mbox{} \\
		& = & [[(0 \to c') \to (0 \to a)']' \to [b \to (c \to a)]']' & \mbox{by (\ref{031114_06})} \\
		& = & [[(0 \to (0 \to c')) \to (0 \to a)']' \to [b \to (c \to a)]']' & \mbox{by (\ref{031114_04})} \\
		& = & [[((0 \to a)' \to (0 \to c')) \to (0 \to a)']' \to [b \to (c \to a)]']' & \mbox{by (\ref{291014_10})} \\
		& = & [[((0 \to a)' \to (0 \to c')) \to (0 \to (0 \to a))']' \to [b \to (c \to a)]']' & \mbox{by (\ref{031114_04})} \\
		& = & [[(0 \to c')' \to (0 \to a)] \to [b \to (c \to a)]']' & \mbox{by (I)} \\
		& = & [[(0 \to c')' \to (0 \to a'')] \to [b \to (c \to a)]']' & \mbox{} \\
		& = & [[0 \to [(0 \to c')'' \to a']'] \to [b \to (c \to a)]']' & \mbox{by (\ref{071114_01})} \\
		& = & [[0 \to [(0 \to c') \to a']'] \to [b \to (c \to a)]']' & \mbox{} \\
		& = & [[0 \to [(c \to 0') \to a']'] \to [b \to (c \to a)]']' & \mbox{by Lemma \ref{general_properties} (\ref{cuasiConmutativeOfImplic})} \\
		& = & [[0 \to [(c \to a'') \to a']'] \to [b \to (c \to a)]']' & \mbox{by (\ref{281014_05})} \\
		& = & [[0 \to [(c \to a) \to a']'] \to [b \to (c \to a)]']' & \mbox{} \\
		& = & [[(c \to a)' \to (0 \to a)] \to [b \to (c \to a)]']' & \mbox{by (\ref{031114_03})} \\
		& = & [(0 \to a) \to b] \to (c \to a) & \mbox{by (I)}
		\end{array}
		$$
		
		\item[(\ref{250315_03})]
		$$
		\begin{array}{lcll}
		a \to ((b \to a) \to b)  & = & a \to [(0 \to a) \to b]  & \mbox{by (\ref{291014_10})} \\
		& = & [(0 \to b) \to a] \to [(0 \to a) \to b] & \mbox{by (\ref{250315_02})} \\
		& = & [(a \to b) \to a] \to [(b \to a) \to b] & \mbox{by (\ref{291014_10})}  \\
		& = & a \to b & \mbox{by (\ref{271114_03})}
		\end{array}
		$$
		
		\item[(\ref{250315_04})]
		$$
		\begin{array}{lcll}
		[(a \to b) \to (b \to c)]'  & = & [[(b \to c)' \to a] \to [b \to (b \to c)]']'' & \mbox{by (I)} \\
		& = & [(b \to c)' \to a] \to [b \to (b \to c)]' & \mbox{} \\
		& = & [(b \to c)' \to a] \to (b \to c)' & \mbox{by (\ref{031114_04})} \\
		& = & [0 \to a] \to (b \to c)' & \mbox{by  (\ref{291014_10})}
		\end{array}
		$$
		
		\item[(\ref{250315_05})]
		$$
		\begin{array}{lcll}
		(a \to b) \to (b \to a) & = & [[(b \to a)' \to a] \to [b \to (b \to a)]']' & \mbox{by (I)} \\
		& = & [[(b \to a)' \to a] \to [b \to a]']' & \mbox{by (\ref{031114_04})} \\
		& = & [[0 \to a] \to [b \to a]']' & \mbox{by (\ref{291014_10})} \\
		& = & (b \to a)'' & \mbox{by (\ref{291014_06})} \\
		& = & b \to a & \mbox{}
		\end{array}
		$$
		
		\item[(\ref{250315_06})]

		
		
		
		$$
		\begin{array}{lcll}
		[[(a \to b) \to c] \to (c' \to a)'] \to [(a \to b) \to c]  & = & [0 \to (c' \to a)'] \to [(a \to b) \to c] & \mbox{by (\ref{291014_10})} \\
		& = & [(c' \to a) \to 0'] \to [(a \to b) \to c] &  \\
		&  & \mbox{by Lemma \ref{general_properties} (\ref{cuasiConmutativeOfImplic})}  &  \\
		& = & [(c' \to a) \to [(a \to b) \to c]'] \to [(a \to b) \to c] & \mbox{by (\ref{281014_05})} \\
		& = & [(c' \to a) \to [(a \to b) \to c]']'' \to [(a \to b) \to c] & \mbox{} \\
		& = & [(a \to (a \to b)) \to c]' \to [(a \to b) \to c]  & \mbox{by (I)} \\
		& = & [(a \to b) \to c]' \to [(a \to b) \to c] & \mbox{by (\ref{031114_04})} \\
		& = & (a \to b) \to c & \\
		&  & \mbox{by Lemma \ref{general_properties_equiv} (\ref{LeftImplicationwithtilde})} &
		\end{array}
		$$
		
		\item[(\ref{250315_07})] 
		\ \\
		$ [[(a \to b) \to c] \to [c' \to (b \to a)]'] \to [(a \to b) \to c] $  \mbox{} \\
		\ \\
		= $\left[[((b \to a) \to (a \to b)) \to c] \to [c' \to (b \to a)]'] \to [((b \to a) \to (a \to b)) \to c\right]  $ 
		\quad \mbox{by (\ref{250315_05})} \\
		\ \\
		 = $((b \to a) \to (a \to b)) \to c $ \quad  \mbox{by (\ref{250315_06})} \\
		 \ \\
		 = $ (a \to b) \to c$   \mbox{by (\ref{250315_05})}. 
		
		 
			
	\end{enumerate}
\end{Proof}

\end{document}